\documentclass[12pt,reqno]{amsart}
\usepackage{setspace}
\usepackage{amsmath, amssymb, amsthm}
\usepackage{amsmath}
\usepackage{amscd,amsthm,amsfonts,amsopn,amssymb,mathrsfs, latexsym}
\usepackage{epsfig,hyperref}
\theoremstyle{plain}
\newtheorem{theorem}{Theorem}

\newtheorem{proposition}[theorem]{Proposition}

\newtheorem{lemma}[theorem]{Lemma}
\theoremstyle{remark}

\allowdisplaybreaks
\def\be{\begin{equation}}
\def\ee{\end{equation}}
\parskip=\medskipamount

\def\ve{\varepsilon}
\def\vp{\varphi}

\def\arrowk{^\to{\kern -6pt\topsmash k}}
\def\arrowK{^{^\to}{\kern -9pt\topsmash K}}
\def\arrowr{^\to{\kern-6pt\topsmash r}}
\def\arrowvp{^\to{\kern -8pt\topsmash\vp}}
\def\arrowf{^{^\to}{\kern -8pt f}}
\def\arrowg{^{^\to}{\kern -8pt g}}
\def\arrowu{^{^\to}a{\kern-8pt u}}
\def\arrowt{^{^\to}{\kern -6pt t}}
\def\arrowe{^{^\to}{\kern -6pt e}}
\def\tk{\tilde{\kern 1 pt\topsmash k}}
\def\barm{\bar{\kern-.2pt\bar m}}
\def\barN{\bar{\kern-1pt\bar N}}
\def\barA{\, \bar{\kern-3pt \bar A}}

\def\iint{\not \kern-4pt\int}
\begin{document}
\title{A MODULAR SZEMER\'EDI-TROTTER THEOREM FOR HYPERBOLAS \\ \
\\ UN TH\'EOR\`EME DE TYPE SZEMEREDI--TROTTER MODULAIRE POUR HYPERBOLES}

\author{Jean Bourgain}
\address{School of Mathematics, Institute for Advanced Study, 1
Einstein Drive, Princeton, NJ 08540.}
\email{bourgain@math.ias.edu}
\thanks{The research was partially supported by NSF grants DMS-0808042 and
DMS-0835373.}

\begin{abstract}
We establish a Szemer\'edi-Trotter type result for hyperbolas in $\mathbb F_p\times \mathbb F_p$.
\end{abstract}
\maketitle

\renewcommand{\abstractname}{R\'esum\'e}
\begin{abstract}
Nous d\'emontrons une version du th\'eor\`eme de Szemer\'edi-Trotter pour des familles d'hyperboles dans $\mathbb F_p\times\mathbb F_p$.
\end{abstract}
\maketitle

\bigskip
\section*{Version fran\c caise abr\'eg\'ee}

Le th\'eor\`eme classique de Szemer\'edi-Trotter donne une estim\'ee sur les incidences d'une famille finie $P$ de points dans le plan et une
famille finie $L$ de droites (ou, plus g\'en\'eralement, de courbes alg\'ebriques de degr\'e born\'e);
dans sa g\'en\'eralit\'e, cette estim\'ee est optimale.
Une version `corps fini', pour $L$ consistant de droites, est obtenue dans \cite {B-K-T}.
Nous nous proposons ici d'\'etablir un r\'esultat de ce type pour certaines familles d'hyperboles dans $\mathbb F_p\times\mathbb F_p$,
d\'efinie par des \'equations
$$
cxy-ax+dy-b=0 \text { o\`u } \begin{pmatrix} a&b\\ c&d\end{pmatrix} \in SL_2(p).
$$
Essentiellement, la condition impos\'ee sur $L$ est que l'ensemble des matrices $\begin{pmatrix} a&b\\c&d\end{pmatrix}$ qui d\'efini $L$ ne
soit pas contenu dans une translat\'ee d'un sous-groupe propre de $SL_2(p)$; l'argument repose en effet sur les r\'esultats de \cite{B-G}
sur l'expansion dans $SL_2(p)$.

\section
{Introduction and statement of the results}

It is shown in \cite{B-K-T} that if $P$ and $L$ are sets of points and lines in $\mathbb P^2(\mathbb F_p)$ with $|P| =|L|=n< p^{2-\ve}$, then
the number of incidences
\be
I(P,L)=|\{(p, \ell)\in P\times L; p\in\ell\}|<cn^{\frac 32 -\delta}
\ee
where $\delta =\delta(\ve)>0$.
An explicit quantitative version of this result appears in \cite{H-R}.
Some of its various applications may be found in \cite {B}.
The following statement provides a result of a similar flavor for hyperbolas.

\begin{proposition}\label{Proposition1}
For all $\ve>0$ and $r>1$, there is a $\delta>0$ such that the following holds.
Let $p$ be a large prime and $A\subset \mathbb F_p$, $S\subset SL_2(p)$ satisfy the conditions
\begin{enumerate}
\item[(2)] $1\ll|A|< p^{1-\ve}$
\item [(3)] $\log|A|< r\log |S|$
\item [(4)] $|S\cap gH|<|S|^{1-\ve}$ for any proper subgroup $H\subset SL_2(p)$ and $g\in SL_2(p)$.
\end{enumerate}

For $g=\begin{pmatrix} a&b\\c&d\end{pmatrix} \in GL_2(p)$, denote $\Gamma_g\subset\mathbb F_p^2$ the curve
$$
cxy-ax+dy-b=0.\eqno{(5)}
$$ 
Then
$$
\{(x, y, g)\in A\times A\times S; (x, y)\in\Gamma_g\}|< |A|^{1-\delta}|S|.\eqno{(6)}
$$
\end{proposition}

Note that the conclusion would be obviously false if we removed assumption (4).

\begin{proposition}\label{Proposition2}

Assume given a polynomial function $\Phi$ on $\mathbb F_p$ taking values in 
{\rm Mat}$_2(\mathbb F_p)$, such that {\rm det}\,$\Phi$ does not vanish
identically and is a quadratic residue (or a quadratic non-residue).
Assume further that $\text{\rm Im\,}\Phi \cap GL_2(p)$ is not contained in a set $\mathbb F^*_p.gH$ for some $g\in SL_2(p)$ and $H\subset SL_2(p)$
a proper subgroup.

Given $\ve>0, r>1$, there is $\delta>0$ such that if $A\subset\mathbb F_p$, $L\subset\mathbb F_p$ satisfy
\begin{enumerate}
\item[(7)] $1\ll |A|< p^{1-\ve}$
\item[(8)] $ \log|A|< r\log L$.
\end{enumerate}

Then
$$
|\{(x, y, t)\in A\times A\times L; (x, y)\in\Gamma_{\Phi(t)}\}|<|A|^{1-\delta}|L|.\eqno{(9)}
$$
\end{proposition}

Applications will be discussed elsewhere.

\section
{Preliminaries}

The main ingredients in the proof are the expansion properties in $SL_2(p)$ obtained in \cite {B-G} and based on \cite{H}.

More specifically, we make use of the so-called `$L^2$-flattening Lemma' that we recall next (in a version slightly more general than stated
in \cite{B-G} but similar proof).

\begin{lemma}\label{Lemma3}

Let $\eta$ be a symmetric probability measure on $SL_2(p)$ and $1\ll K< p^{\frac 1{10}}$, such that
\begin{enumerate}
\item[(10)] $\eta(gH)<K^{-1} $ for any proper subgroup $H\subset SL_2(p), g\in SL_2(p)$
\item[(11)] $\Vert \eta\Vert_2> Kp^{-3/2}$.
\end{enumerate}
Then
$$
\Vert\eta*\eta\Vert_2 < K^{-c} \Vert\eta\Vert_2\eqno{(12)}
$$
with $c>0$ an absolute constant.
\end{lemma}

Here $\Vert\eta\Vert_2 =\Big[\sum_{g\in SL_2(p)}\eta(g)^2\Big]^{\frac 12}$.
Denoting $\eta^{(\ell)}$ the $\ell$-fold convolution, iteration of Lemma 3 implies that for $\eta$ satisfying (10)
$$
\Vert\eta^{(2^\ell)}\Vert_2 \leq K^{-c\ell}\Vert\eta\Vert_2 +Kp^{-3/2}\eqno{(13)}
$$
and hence
$$
\Vert\eta^{(2^{\ell+1)}} \Vert_\infty \leq \Vert\eta^{(2^\ell)}\Vert^2_2 \leq 2K^{-2c\ell}\Vert\eta\Vert^2_2 +2K^2 p^{-3}.
\eqno{(14)}
$$
Recall also that if $K>p^\gamma$, (14) combined with an argument due to \cite{S-X} based on Frobenius multiplicity, implies that
$$
\Vert\eta^{(\ell)} \Vert_\infty< 2p^{-3}\eqno{(15)}
$$
for some $\ell =\ell(\gamma)$.
See \cite{B-G} for details.

\section
{Sketch of the proof of Proposition \ref{Proposition1}}

Using the action $\tau$ of $SL_2(p)$ on $\mathbb P^1(\mathbb F_p)$, rewrite equation (5) (since we may assume $cx+d\not=0$) as
$$
y=\frac {ax+b}{cx+d}=\tau_g(x).
$$
The left hand side of (6) equals
$$
\sum_{g\in S} |A\cap \tau_{g^{-1}}(A)|=|S|.\langle1_A, \sum_g(1_A\circ \tau_g)\mu(g)\rangle
$$
with
$$
\mu= \frac 1{|S|} \sum_{g\in S} \delta _g
$$
and $\langle \, , \, \rangle$ referring to the inner product on $L^2(\mathbb F_p)$.

Next, applying the Cauchy-Schwarz inequality, write
$$
\begin{aligned}
\langle 1_A, &\sum(1_A\circ\tau_g)\mu(g)\rangle \leq |A|^{\frac 12} \Big\Vert\sum_g(1_A\circ \tau_g)\mu(g)\Big\Vert_2\\
&=|A|^{\frac 12}\Big[\sum \langle 1_A\circ\tau_g, 1_A\rangle (\mu*\mu^{-1})(g)\Big]^{\frac 12}\\
&\leq |A|^{\frac 34} \Big\Vert \sum(1_A\circ \tau_g)\nu(g)\Big\Vert_2^{\frac 12}
\end{aligned}
$$
where $\nu=\mu*\mu^{-1}$ is symmetric.
Iteration gives for any $\ell\in\mathbb Z_+$
$$
\leq |A|^{1-2^{-\ell-1}}\Big\Vert\sum(1_A\circ \tau_g)\nu^{(2^{\ell-1})}(g)\Big\Vert_2^{2^{-\ell}}.\eqno{(16)}
$$
Note that, by our assumption (4), if $H\subset SL_2(p)$ is a proper subgroup and $g\in SL_2(p)$,
$$
\nu(gH)<|S|^{-\ve}.\eqno{(17)}
$$
It follows from (16) that if (6) fails, then
$$
|A|^{\frac 12-2^\ell\delta}<\Big\Vert \sum(1_A\circ \tau_g)\nu^{(2^{\ell-1})} (g)\Big\Vert_2\eqno{(18)}
$$
and hence
$$
|A|^{1-2^{\ell+1}\delta} <\sum_{x\in \mathbb F_p, g\in SL_2(p)} 1_A(\tau_gx)1_A(x)\nu^{(2^\ell)}(g).\eqno{(18)}
$$

We distinguish two cases.

If $|A|>p^{\frac 1{10}}$, then, by (3), $|S|>p^{\frac 1{10r}}$ and, applying (15) with \hfill\break
$K=p^{\frac\ve{10r}}$ gives
$\vert\nu^{(2^\ell)}\Vert_\infty < 2p^{-3}$ for some $\ell=\ell(\ve, r)$.
Hence, from (18), $|A|^{1-2^{\ell+1}\delta}<2p^{-1}|A|^2$, contradicting (2) for $\delta$ small enough.

If $|A|\leq p^{\frac 1{10}}$, denote $\nu_1=\nu^{(2^\ell)}$ and write using (18) and H\"older
$$
\sum_g\Big[\sum_{x_1, x_2, x_3\in A} 1_A(\tau_gx_1)1_A(\tau_gx_2) 1_A(\tau_gx_3)\Big] \nu_1(g)> |A|^{3-32^{\ell+1}\delta}.\eqno{(19)} 
$$
Assuming $2^\ell\delta <\frac 1{10}$, it follows from (19) that there are distinct elements $x_1, x_2, x_3\in A$ such that
$$
\nu_1[g\in SL_2(p); \tau_g x_i \in A \text { for } i=1, 2, 3]>\frac 12|A|^{-32^{\ell+1}\delta}.\eqno{(20)}
$$

Since the equations $\tau_g x_i = y_i(i= 1, 2, 3)$ determine $g$ up to bounded multiplicity, (20) implies $\Vert\nu_1\Vert_\infty
|A|^3 > c|A|^{-32^{\ell+1}\delta}$.
Since by (14), applied with $K=|A|^{\frac \ve r}, \Vert\nu_1\Vert_\infty\leq 2|A|^{-2c(\ell-1)\frac\ve r} +p^{-2}$, it follows
$$
|A|^{3-2c(\ell-1)\frac \ve r}+ |A|^3 p^{-2}> c|A|^{-32^{\ell+1}\delta}.\eqno{(21)}
$$
Taking $\ell =\ell(\ve, r)$ appropriately and $\delta$ small enough gives again a contradiction.

\section
{Proof of Proposition \ref{Proposition2}}

We may assume that $\det \Phi (t)\in\mathbb F_p$ is a quadratic residue (otherwise replace $\Phi$ by $\Phi(t_0)\Phi$ with $\det \Phi(t_0)\not=0$).
We can also assume that $\det\Phi(t)\not=0$ for $t\in L$.
Apply Proposition \ref{Proposition1} to
$$
S=\{\sigma (t)\Phi(t); t\in L\}\subset SL_2(p)\eqno{(22)}
$$
where $\sigma(t)\in\mathbb F_p^*$ is chosen such that $\sigma(t)^2 \det\Phi(t)=1$.
It remains to verify condition (4).  From our assumption on $\Phi$,
$$
\mathbb F_p^*.\langle \Phi(t_1)^{-1} \Phi(t_2); \det \Phi (t_1) \not= 0, \det \Phi (t_2) \not= 0\rangle =GL_2(p).\eqno{(23)}
$$

Using the bicommutator characterization of large proper subgroups of $SL_2(p)$, (23) implies that also
$$
\mathbb F_p^* \langle \Phi(t_1)^{-1}\Phi(t_2); t_1, t_2 \in L_1\rangle = GL_2(p)\eqno{(24)}
$$
for any sufficiently large subsets $L_1$ of $L$. Hence
$$
\{\sigma(t)\Phi(t); t\in L_1\}
$$
is not contained in a coset of a proper $SL_2(p)$-subgroup.

\end{document}